\documentclass[11pt,oneside,a4paper]{article}

\usepackage{lmodern}
\usepackage{babel}
\usepackage[T1]{fontenc} 
\usepackage[utf8]{inputenc} 
\usepackage{verbatim}  
\usepackage[a4paper]{geometry} 
\usepackage{amsmath} 
\usepackage{graphicx} 
\graphicspath{{images/}} 
\usepackage{systeme}
\usepackage{amssymb}
\usepackage{amsfonts}
\usepackage{csquotes}
\usepackage{hyperref}

\usepackage{amsmath,amssymb}
\usepackage{amsthm}

\newtheorem{theorem}{Theorem}[section]
\newtheorem{conjecture}{Conjecture}[section]

\newtheorem{lemma}{Lemma}[section]
\newtheorem{corollary}{Corollary}[section]

\newtheorem{remark}{Remark}[section]
\newtheorem{remarks}{Remark}[section]
\newtheorem{examples}{Examples}[section]

\usepackage{amssymb}
\usepackage{amsmath}
\usepackage{amstext}

\usepackage{comment}

\newcommand*\modd[3]{#1\equiv#2\;( {\rm mod}\;  #3)}

\title{On a general Syracuse problem with conjectures}

\author{A. Bouhamidi\footnote{L.M.P.A, Universit\'e du
		Littoral, 50 rue F. Buisson BP699, F-62228 Calais-Cedex, France.  \href{mailto:abderrahman.bouhamidi@univ-littoral.fr}{abderrahman.bouhamidi@univ-littoral.fr}}}

\date{August 30,  2021}

\begin{document}
\maketitle

\begin{abstract}
	In this paper, we study a  general Syracuse problem.  We give some  necessary conditions  concerning the existence of  eventual non trivial cycles.   Some properties based on  linear logarithmic forms are established. New general conjectures are given. To illustrate  the behavior of  such a problem,  some  particular examples are presented.
\end{abstract}



\section{Introduction}
 Let $\mathbb{N}=\{1,2,3,\ldots\}$ be the set of integers $\geq 1$.  The Syracuse problem  also known as the Collatz, Kakutani or  3x+1 problem, concerns the sequence of positive integers generated by the  iterations associated to the following map $S:\mathbb{N}:\longrightarrow \mathbb{N}$ defined by
 \begin{equation}\label{CollatzMap}
 	S(n)=\left\{\begin{array}{lcc}
 		n/2&\mathrm{if}& n\;\mathrm{is\; even}\\ 
 		\dfrac{3n+1}{2}&\mathrm{if}& n\;\mathrm{is\; odd}.\\
 	\end{array}\right.
 \end{equation}
For every $k\in\mathbb{N}_0=\mathbb{N}\cup\{0\}$, the notation $S^{(k)}$ stands for the $k$th iterate of the map $S$.  The well-known conjecture of the $3x+1$ problem (or Collatz's conjecture ) states that, for  any initial value $n\in\mathbb{N}$, there is a positive integer $k\in\mathbb{N}_0$ such that $S^{(k)}(n)=1$. A second conjecture in the Syracuse problem states  that the trivial cycle $(1\rightarrow 2\rightarrow1)$ is the unique cycle in the graph associated to the map $S$. The $3x+1$ problem has been explored for more than $80$ years by many authors.
An extensive overviews on the  literature concerning the Syracuse problem has been given by Lagarias \cite{Lagarias1985,Lagarias1990,Lagarias2011}.  The conjecture has been verified experimentally with a computer \cite{Oliviera2010},  for $n\leq 5\times 2^{60}\simeq 5\, 764\, 607\, 523\, 034\, 234\, 880$.  In  \cite{Eliahou1993},  Eliahou  showed that any nontrivial cyclic  under the operator $S$ must contains at least about $1708791$ elements provided that $\min(\Omega)>2^{40}$.  The latest calculations
based on the computation of  $\min(\Omega)>87 \times 2^{60}$ increases minimum cycle length to $10439860591$. Recently, Tao in \cite{Tao2020}  proved that Collatz conjecture is almost true for almost all numbers.
Many generalizations of the Syracuse problem have been given in the literature, see \cite{Lagarias2011} for more details.

In this paper, we study a generalization to the Syracuse problem. The goal is to shed light on the sequence related to the $3x+1$ problem for a better understanding of the phenomenon.  The general problem $"qn+r"$ was first briefly introduced by Crandall in  \cite{Crandall1978}. In this paper \cite{Crandall1978}, Crandall mentioned some experimental  results and gave some discussions on a such general problem.    In our present paper, we will study  such a problem in more details and we will give some new results.

 The outline of this paper is as follows. In Section 2, we introduce the unified   generalization of the Syracuse problem. We will give some notations, definitions and results  relative to the general operator. In Section 3, we will give some results on the  constraints on eventual non trivial cycle. This section is in fact a generalization of the  work of Eliahou  \cite{Eliahou1993}. Section 4, we study a few  examples with some  results and conjectures.

\section{The Syracuse generalized problem}\label{SectionGeneralization}
We consider the map  $T:\mathbb{N}\longrightarrow \mathbb{N} $ given  by
\begin{equation}\label{mapT}
	T(n)=\left\{\begin{array}{lc}
	n/2&\mathrm{if }\; $n$\; \mathrm{ is\; even},\\
	(a n+ b )/2&\mathrm{if} \;$n$\;\mathrm{ is\;odd},\\
	\end{array}\right.
\end{equation}
where $a$ and $b$ are both odd integers such that $a>-b$.  Then, the number $an+b$ is an even positive integer for any odd number $n$ and one can iterate the $T$ function any number of times.

The Syracuse  problem may be  viewed as  a direct graph whose vertices are the positive integers and whose edges are the connection from $n$ to $T(n)$. Following the terminology of Lagarias \cite{Lagarias1985},  such a graph will be called the Collatz graph. Given $n\in \mathbb{N}$, the  trajectory of $n$ is the set $\Gamma(n)=\{n,T(n),T^{(2)},\ldots\}$ of iterates.  A cycle
of the graph associated to  the map $T$ having $k$ vertices   is a trajectory $\Omega$  for which $T^{(k)}(x)=x$ for all $x\in \Omega$.  The cardinal $\#\Omega=Card(\Omega)=k$, also called the length, of $\Omega$ is in fact the smallest integer such that $T^{(k)}(x)=x$ for all $x\in \Omega$.  We will denote by  
$\Omega(\omega)=\{\omega,T(\omega),T^{(2)}(\omega),\ldots, T^{(k-1)}(\omega)\}$
the  cycle associated to $T$  of length  $k$
where $\omega$ is the smallest element of the cycle. We will also use the notation 
$$\Omega(\omega)=(\omega \rightarrow T(\omega)\rightarrow T^{(2)}(\omega)\rightarrow \ldots \rightarrow T^{(k-1)}(\omega)\rightarrow \omega),$$
to denote the cycle $\Omega(\omega)$ having $k$ elements. For each choice of the parameters $a$ and $b$ the notation $\mathcal{C}(a,b)$ stands for the set of all possible cycles associated to $T$. 

There are two possible situations for a  trajectory $\{T^{(k)}(n)\}_{k\in\mathbb{N}_0}$, for $n\in\mathbb{N}$:
\begin{enumerate}
	\item[(i)]  Convergent trajectory to one or other  specific  cycle $\Omega(\omega)$, assuming that   such a cycle exists.
	\item[(ii)] Divergent trajectory: $\displaystyle \lim_{k\rightarrow \infty}T^{(k)}(n)=\infty$.
\end{enumerate}

Assuming that the cycle $\Omega(\omega)$ exists, we denote by $G(\omega)$ and $G(\infty)$ the subsets of $\mathbb{N}$  given by
$$G(\omega)=\{n\in\mathbb{N}\;:\; \exists k\in\mathbb{N},\quad T^{(k)}(n) =\omega\},$$
and 
$$G(\infty)=\{n\in\mathbb{N}\;:\; \lim_{k\rightarrow \infty} T^{(k)}(n)=\infty\}.$$
We observe that for all $\omega,\omega'\in\mathbb{N}$ such that the  cycles $\Omega(\omega)$ and  $\Omega(\omega')$ exist, if  $\omega\not= \omega'$, then  $G(\omega)\cap G(\omega')=\emptyset$ and 
$G(\infty)\cap G(\omega)=\emptyset$. The sets are pairwise  disjoint and the set $G(\infty)$ may be empty. We have $\Omega(\omega)\subset G(\omega)$.

\begin{conjecture}\label{GeneralConjecture}For any parameters $a$ and $b$ such that $a+b=2^{\nu_0}$, where $\nu_0\geq 1$,  the set  $\mathcal{C}(a,b)$ has  finite cardinality: 
	$1\leq q^*=\#\mathcal{C}(a,b)<\infty$.
Then, it follows that  there is  $q^*$ distinct integers $\omega_1^*,\ldots,\omega_{q^*}^*$ such that:
	$\mathcal{C}(a,b)=\left\{\Omega(\omega_1^*),\ldots,\Omega(\omega_{q^*}^*)\right\}$,
	and
	$$\mathbb{N}=G(\omega_1^*)\cup\cdots\cup G(\omega_{q^*}^*)\cup  G(\infty).$$
	The set $G(\infty)$ may be empty or not empty.
\end{conjecture}

\textbf{Questions:} Suppose that the condition $a+b=2^{\nu_0}$ holds, then one can ask  the following questions: 
\begin{enumerate}
	\item  Under which conditions  on  the parameters $a$ and $b$  the set $ G(\infty)$ is empty or not empty? 
	In this case where the set $ G(\infty)$ is empty, there is no divergent trajectory.
	\item Is it possible to find an  estimate for  the cardinality  $\#\mathcal{C}(a,b)$? Then, it will be possible to confirm that the set  $\mathcal{C}(a,b)$  is finite. 
	\item  Under which conditions  on  the parameters $a$ and $b$, the set $ \mathcal{C}(a,b)$ is of cardinality $ \#\mathcal{C}(a,b)=1$ ? In this case, we will see that the unique cycle in $\mathcal{C}(a,b)$  is the trivial one $\Omega(1)$.  If in addition $G(\infty)=\emptyset$, then for any integer $n$ there exists an integer $k$ such that $T^{(k)}(n)=1$ and we will  have $\mathbb{N}=G(1)$.

\end{enumerate}
To  find  answers to the previous  questions in a general situation, is a difficult task. We will study in the following sections the problem in some particular situations. 
The classical conjecture of Collatz may be reformulated as following.
\begin{conjecture}\label{Collatzconjecture} (Collatz conjecture) For  the parameters $a=3$ and $b=1$  the set  $\mathcal{C}(3,1)$  contains only one cycle (the trivial one): 
	$\mathcal{C}(3,1)=\left\{\Omega(1)\right\}$,
	with $\Omega(1)=(1\rightarrow 2 \rightarrow 1)$ and we  have: 
	$\mathbb{N}=G(1)$ and $G(\infty)=\emptyset$.
	
\end{conjecture}

\begin{remark}
	If we assume that  the set  $\mathcal{C}(a,b)$ is of  finite cardinality. Then,  the sets 
	$G(\omega^*_1), \ldots,G(\omega^*_q),G(\infty)$ are the classes for the
binary equivalence relation $\sim_T$ defined on $\mathbb{N}$ as following: For $n,n'\in\mathbb{N}$
	$$n \sim_T n'\Longleftrightarrow \exists\, \omega\in \mathbb{N}\cup\{\infty\}:\;\; n,n'\in G(\omega).$$
The quotient of $\mathbb{N}$ by the relation $\sim_T$ is 
	$$\mathbb{N}/\sim_T=\Bigl\{ G(\omega^*_1), \ldots,G(\omega^*_q),G(\infty)\Bigr\},$$
	if  $G(\infty)\not=\emptyset$ and 
	$$\mathbb{N}/\sim_T=\Bigl\{ G(\omega^*_1), \ldots,G(\omega^*_q)\Bigr\},$$
	if   $G(\infty)=\emptyset$. Furthermore,
$$\mathbb{N}=G(\omega_1^*)\cup\cdots\cup G(\omega_q^*)\cup  G(\infty).$$
\end{remark}

In the following theorem, we will see that if $a$ and $b$ satisfy the  condition of Conjecture~\ref{GeneralConjecture}, then    the set  $\mathcal{C}(a,b)$ is not empty.  

\begin{theorem}\label{SetOfCycles}For any parameters $a$ and $b$ such that $a+b=2^{\nu_0}$, where $\nu_0\geq 1$,  the set  $\mathcal{C}(a,b)$ is not empty. It contains at least the  trivial cycle of length $\nu_0$:
	$$\Omega(1)=(1 \rightarrow  2^{\nu_0-1}\rightarrow 2^{\nu_0-2}\rightarrow \ldots \rightarrow 2 \rightarrow 1).$$
	Furthermore, if in addition  $a=2^{\nu_1}-\delta$, where $\nu_1\geq 1$ for $\delta=\pm1$ with  $\delta b>0$, then the set $\mathcal{C}(a,b)$ contains in addition the following  second  trivial cycle of length $\nu_1$: 
	$$\Omega(\delta b)=(\delta b \rightarrow  \delta 2^{\nu_1-1} b\rightarrow \delta2^{\nu_1-2}b\rightarrow \ldots \rightarrow 2\delta  b \rightarrow \delta b).$$ 
\end{theorem}
\begin{proof}
	It is trivial to check that for all odd integers  $a$ and $b$  such that $a+b=2^{\nu_0}$, we have $T^{(\nu_0)}(1)=~1$. Indeed, we have
	$T(1)=(a+b)/2=2^{\nu_0-1}$ and then $T^{(\nu_0)}(1)=~1$. 
	Furthermore, if in addition  $a=2^{\nu_1}-\delta$  for  $\delta=\pm1$ with  $\delta b>0$, then we have $T^{(\nu_1)}(\delta b)=~\delta b$. Indeed, $T(\delta b)=(\delta ab+b)/2=b (\delta a+1)/2$. But as $\delta a=\delta 2^{\nu_1}-1$. It follows that $T(\delta b)=(\delta b) 2^{\nu_1-1}$ and then $T^{(\nu_1)}(\delta b)=\delta b$.
\end{proof}

\begin{remark}
	According to the previous Theorem, we  may set $\omega_1^*=1$. Furthermore,  in the case where $\delta b>0$, we have  $\omega_2^*=\delta b$ and $ \#\mathcal{C}(a,b)\geq 2$.
\end{remark}

\begin{examples}
	\begin{enumerate}
		\item For  the classical Syracuse problem, 
		we have $a=3=2^{\nu_1}+\delta$ and  $b=1=2^{\nu_0}-a$, where
 $(\nu_0,\nu_1,\delta)=(1,2,1)$. In this case, $\mathcal{C}(3,1)$ has at least the  trivial  cycle  of length $\nu_1=2$: $\Omega(1)=(1\rightarrow 2\rightarrow1).$
		
		We observe that $(a,b)=(3,1)$ may be also written in the form 
		 $a=3=2^{\nu_1}-\delta$ and $b=2^{\nu_0}-a=1$ with $(\nu_0,\nu_1,\delta)=(2,2,1)$ and we have $\delta b>0$.  So according to the previous theorem, $\mathcal{C}(3,1)$  must contains two  trivial cycles $\Omega(1)$ of length $\nu_1$ and $\Omega( \delta b)$ of length $\nu_0$.  Here of course the two trivial cycles coincide  $(\delta b=1)$.

		\item For  $a=7$ and $b=9$, 	the corresponding map $T$ is given by
		$$
		T(n)=\left\{\begin{array}{lc}
			n/2&\mathrm{if }\; $n$\; \mathrm{ is\; even},\\
			(7 n+9 )/2&\mathrm{if} \;$n$\;\mathrm{ is\;odd}.\\
		\end{array}\right.
		$$
		We have $a=7=2^{\nu_1}-\delta$ and $b=9=2^{\nu_0}-a$ with $(\nu_0,\nu_1,\delta)=(4,3,1)$ and $\delta b>0$.
		According to the previous theorem, $\mathcal{C}(7,9)$  contains at least   the two trivial cycles the first one $\Omega(1)$ of length $\nu_0=4$ and the second one  $\Omega(\delta b)$ of length $\nu_1=3$. We have: 
		$\Omega(1)=(1\rightarrow 8\rightarrow 4\rightarrow 2\rightarrow1)$ and $\Omega(\delta b)=(9\rightarrow 36\rightarrow 18\rightarrow9)$.
		
		\item For  $a=9$ and $b=-7$, 	the corresponding map $T$ is given by
				$$
		T(n)=\left\{\begin{array}{lc}
			n/2&\mathrm{if }\; $n$\; \mathrm{ is\; even},\\
			(9 n-7 )/2&\mathrm{if} \;$n$\;\mathrm{ is\;odd}.\\
		\end{array}\right.
		$$
		We have $a=9=2^{\nu_1}-\delta$ and $b=-7=2^{\nu_0}-a$ with $(\nu_0,\nu_1,\delta)=(3,1,-1)$ and $\delta b>0$. According to the previous theorem, $\mathcal{C}(9,-7)$  contains at least   the two trivial cycles the first one $\Omega(1)$ of length $\nu_1=1$ and the second one  $\Omega(\delta b)$ of length $\nu_0=3$. We have: 
		$\Omega(1)=(1\rightarrow1)$ and $\Omega(\delta b)=(7\rightarrow 28\rightarrow 14\rightarrow7)$.

\item For  $a=3$ and $b=5$, 	the corresponding map $T$ is given by
$$
T(n)=\left\{\begin{array}{lc}
	n/2&\mathrm{if }\; $n$\; \mathrm{ is\; even},\\
	(3 n+5 )/2&\mathrm{if} \;$n$\;\mathrm{ is\;odd}.\\
\end{array}\right.
$$
We have $a=3=2^{\nu_1}-\delta$ and $b=5=2^{\nu_0}-a$ with $(\nu_0,\nu_1,\delta)=(3,2,1)$ and $\delta b>0$.
According to the previous theorem, $\mathcal{C}(3,5)$  contains at least   the two trivial cycles. The first one $\Omega(1)$ of length $\nu_0=3$  and the second one  $\Omega(\delta b)$ of length $\nu_1=2$. We have: 
		$\Omega(1)=(1\rightarrow 4\rightarrow 2\rightarrow1)$ and $\Omega(\delta b)=(5\rightarrow 10\rightarrow 5)$.
		But, in this case, $\mathcal{C}(3,5)$ has in addition two non trivial cycles of lengths $5$, which are:
		
		$\Omega(19)=(19\rightarrow 31\rightarrow 49\rightarrow 76\rightarrow38\rightarrow19)$,
		
		$\Omega(23)=(23\rightarrow 37\rightarrow 58\rightarrow 29\rightarrow46\rightarrow23)$, and
		two other non trivial cycles of lengths $27$, which are:
		
		$\Omega(187)=(187\rightarrow283\rightarrow427\rightarrow643\rightarrow967\rightarrow1453\rightarrow2182\rightarrow1091\rightarrow1639\rightarrow2461\rightarrow3694\rightarrow1847\rightarrow2773\rightarrow4162\rightarrow2081\rightarrow3124\rightarrow1562\rightarrow781\rightarrow1174\rightarrow587\rightarrow883\rightarrow1327\rightarrow1993\rightarrow2992\rightarrow1496\rightarrow748\rightarrow374\rightarrow187)$ and 
		
		$\Omega(347)=(347\rightarrow523\rightarrow787\rightarrow1183\rightarrow1777\rightarrow2668\rightarrow1334\rightarrow667\rightarrow1003\rightarrow1507\rightarrow2263\rightarrow3397\rightarrow5098\rightarrow2549\rightarrow3826\rightarrow1913\rightarrow2872\rightarrow1436\rightarrow718\rightarrow359\rightarrow541\rightarrow814\rightarrow407\rightarrow613\rightarrow922\rightarrow461\rightarrow694\rightarrow347)$.
	\end{enumerate}

\end{examples}

The arguments for the following theorem are similar to those   from  \cite{Crandall1978}.

\begin{theorem}
For $b>1$,  $\exists n\in\mathbb{N}$ such that  $\forall k\in \mathbb{N}_0 $, we have $T^{(k)}(n)\not=1$. 
\end{theorem}
\begin{proof}
	Let $b>1$ be an odd integer and choose any integer $n$ such that $\modd{n}{0}{b}$. Then
	$2T(n)=\modd{b(aq+1)}{0}{b}$, where
	$n=qb$. It follows that  $\modd{2^\ell T(n)}{0}{b}$ for any power $\ell\geq 1$. As $b$ is odd,  by Gauss Theorem , we have also
	$\modd{T(n)}{0}{b}$. It follows, that all iterates  $T^{(k)}(n)$ of $n$ are  multiple of $b$. Hence, for 
	$\modd{n}{0}{b}$ and $b>1$, there is no integer $k$ such that $T^{(k)}(n)=1$.
\end{proof}

\begin{corollary}
If  $b> 1$ and   $\#\mathcal{C}(a,b)=1$,  then $G(\infty)\not=\emptyset$, i.e., the map $T$ has a divergent trajectory.
\end{corollary}
\begin{proof}
	Its an immediat consequence of the previous theorem.
\end{proof}

\section{Bounds for cycle lengths}

Let $\xi$ be an irrational real number. We say that $\mu= \mu(\xi)$  is an effective irrationality measure  of $\xi$  if for all $\varepsilon>0$,  there exists an integer   $q_0(\varepsilon)>0$  (effectively computable) such that   
$$\forall (p,q) \in \mathbb{Z}\times \mathbb{N},\quad   q>q_0(\varepsilon) \Longrightarrow
\Bigl|\dfrac{p}{q}-\xi\Bigr|>\dfrac{1}{q^{\mu+\varepsilon}},$$
see \cite{Rhin1987} for more details.   The following theorem  is a  consequence of Baker’s theory of linear forms in logarithms (see e.g.  \cite{Bugeaud2015, Waldschmidt2000} and the references therein). By definition, two positive rational numbers are multiplicatively independent if the quotient of their logarithms is irrational.  In the following, we denote by $\xi$ 
the  real number $\xi=\dfrac{\log(a)}{\log(2)}$. We note that  $\xi$ is a transcendental real number.  Indeed, suppose that  $\xi$ is a rational number, then there exist a pair of integers $(p,q)\in\mathbb{N}^2$ such that $a^q=2^q$, which is impossible since $a$ is an odd integer number.
Now, as $\xi$ is irrational real number and $2^\xi=a$ is algebraic number, then according to the Gelfond-Schneider theorem, we get that $\xi$ is transcendental real number.
\begin{theorem}
	Let $a_1$, $a_2$, $b_1$, $b_2$ be positive integers with $a_1>a_2$ and $b_1>b_2$. Assume that $a_1/a_2$ and $b_1/b_2$ are multiplicatively independent. There exists an absolute, effectively computable, constant $C$
	such that
	$$\mu_{eff}\Bigl(
	\dfrac{\log(a_1/a_2)}{log(b_1/b_2) }\Bigr)\leq C(\log a_1 )(\log b_1 ).$$
\end{theorem}

In our situation, we set   $a_1=a>a_2=1$ and $b_1=2>b_2=1$. 
According to the previous theorem, we have
	$\mu_{eff}(\xi)\leq C(\log a )(\log 2 )$.
The  convergents of the continued fraction of the irrational number $\xi$ imply that $\mu_{eff}(\xi) \geq  2$. Thus, the effective irrationality measure  of $\xi$ is a finite number:
$$2\leq \mu_{eff}(\xi)\leq C(\log a )(\log 2 ).$$
In the following, we set $\mu=\mu_{eff}(\xi)$.

For $\varepsilon=1$,  there exists an integer $\widehat{q}_0>0$ such that
\begin{equation}\label{muequation}
\forall (p,q) \in \mathbb{Z}\times \mathbb{N},\quad   q>\widehat{q}_0 \Longrightarrow
\Bigl|\dfrac{p}{q}-\xi\Bigr|>\dfrac{1}{q^{\mu+1}}\Longrightarrow \Bigl|p-q\xi\Bigr|>\dfrac{1}{q^{\mu}}.
\end{equation}
Let us now recall some properties from the continued fractions which may be found in the literature of the theory of rational approximation. 
There is a unique representation of $\xi$ as an infinite simple continued fraction:
$$\xi=[a_0,a_1,a_2,a_3,\ldots]$$
The integer numbers $a_n$ are the partial quotients, the rational numbers 
$$\dfrac{p_n}{q_n} = [a_0,a_1,\ldots,a_n]$$
are the convergents, where $gcd(p_n,q_n)=1$. The numbers
$x_n= [a_n,a_{n+1},\ldots]$
are the complete quotients.  From these definitions we deduce, for $n \geq 0$,
$$\xi= [a_0,a_1,\ldots,a_n,x_{n+1}]=\dfrac{x_{n+1}p_n+p_{n-1}}{x_{n+1}q_n+q_{n-1}}.$$
 The integers $p_n,q_n$ are obtained recursively as follows:
 $$\begin{array}{rcl}
p_{n}&=&a_np_{n-1}+p_{n-2},\\
q_{n}&=&a_nq_{n-1}+q_{n-2},\\
 \end{array}
$$
 with the initial values $p_{-2}=0$, $p_{-1}=1$ and  $q_{-2}=1$, $1_{-1}=0$.
The sequences $(p_n)$ and $(q_n)$ are strictly increasing unbounded sequences.  It is well known that 
 the fractions $p_n/q_n$ satisfy the following properties:

\begin{lemma}\label{Lemma1ContinuedFractions} For any pair of integers $(p,q)\in\mathbb{N}^2$ with $1<q<q_n$,
$$
		\Bigl|p_n-q_n\xi\Bigr|<\Bigl|p-q\xi\Bigr|.
$$
\end{lemma}

\begin{lemma}\label{Lemma2ContinuedFractions} For $n\geq 0$,
	$$ \dfrac{1}{q_n+q_{n+1}}<
	\Bigl|{p_n}-{q_n}\xi\Bigr|.
	$$
\end{lemma}

  Let $\Omega$ be a cycle of $T$ and let $\Omega_0$ and $\Omega_1$ denote the subsets of $\Omega$ of its even and odd elements, respectively. Let   $L=\#\Omega_0$ and  $K=\#\Omega_1$ be the cardinals of $\Omega_0$ and $\Omega_1$, respectively. We have the following theorem which may be seen as an extension to the one obtained by Eliahou \cite{Eliahou1993}.

\begin{lemma}\label{theoremEncadrement(K+L)-K}
	We assume $b\geq 1$. Then,
$K$ and $L$ satisfy  the following inequality
\begin{equation}\label{log(a)log(2)}
 0	<(K+L)-K\xi\leq 
\dfrac{b\, K}{a\log(2)\min(\Omega_1)},
\end{equation}
wish implies that
\begin{equation}\label{log(a)log(2)}
	0	<(K+L)-K\xi\leq \dfrac{b\;\; \#\Omega}{a\log(2)\min(\Omega)} .
\end{equation}
\end{lemma}
\begin{proof}
We have 
$$\Bigl(\prod_{n\in \Omega}n\Bigr)= \Bigl(\prod_{n\in \Omega}T(n)\Bigr)= \Bigl(\prod_{n\in \Omega_0}T(n)\Bigr) \times \Bigl(\prod_{n\in \Omega_1}T(n)\Bigr).$$
Dividing the last equality by the member in the left hand side, we get
$$1<\dfrac{2^{K+L}}{a^K}=\prod_{n\in \Omega_1} \Bigl(1+\dfrac{b}{an}\Bigr)\leq \Bigl(1+\dfrac{b}{a\min\Omega_1}\Bigr)^K .$$
Applying the $\log$ and the fact that $\log(1+x)<x$ for $x>0$. We obtain the required result.
\end{proof}

 The  following theorem  is an extension to a result obtained by  Crandall \cite{Crandall1978}.
\begin{theorem}\label{theorem1_K>=}
	Let us set   $c_0 =\dfrac{a\log(2)}{b}$ with $b\geq 1$. Then,  for all $n\geq 1$, 
	$$ K\geq\min\Bigl(q_n,\dfrac{c_0 \min(\Omega_1)}{q_n+q_{n+1}}\Bigr),$$
	which implies that
		$$\#\Omega\geq\min\Bigl(q_n,\dfrac{c_0 \min(\Omega)}{q_n+q_{n+1}}\Bigr),$$
\end{theorem}
\begin{proof}
	Let $n>1$, if $K>q_n$, then the  required inequality  is obvious. Assume that $K\leq q_n$, it follows from  Lemma~\ref{Lemma1ContinuedFractions}  that 
	$|(K+L)-K\xi|>|p_n-q_n\xi|$. By Lemme~\ref{Lemma2ContinuedFractions}, we get
	$$|(K+L)-K\xi|>\dfrac{1}{q_n+q_{n+1}}.$$
	And now, from Lemma~\ref{theoremEncadrement(K+L)-K}, we get
	$$\dfrac{b\, K}{a\log(2)\min(\Omega_1)}>\dfrac{1}{q_n+q_{n+1}}.$$
	Which prove the requried inequalities.
\end{proof}

\begin{remarks}
	\begin{enumerate}
		\item 	We observe that for the classical case with  $a=3$ and $b=1$, we have $c_0\simeq 2.07944\ldots\geq 2$ which the coefficient appearing in \cite{Crandall1978}.
		\item  As an application of the previous theorem. For $a=3$ and $b=1$ we know that the Collatz conjecture  has been verified experimentally with a computer \cite{Oliviera2010},  until  $n\leq 
		N_0=5\times 2^{60}\simeq 5\, 764\, 607\, 523\, 034\, 234\, 880$. With
		$q_{19}=397\,573\,379 $ and 	$q_{20}=6\,189\,245\,291$, the cardinal  of an hypothetical nontrivial cycle $\Omega$  must  satisfy
		$$\#\Omega \geq 3.63974\times 10^8\simeq 363\,974\,000.$$ Similar result was obtained by Eliahou in \cite{Eliahou1993} by using Farey fractions and also continued fractions.
	\end{enumerate}
\end{remarks}

\begin{theorem}\label{theorem2_K>=}
		Let us set   $c_0 =\dfrac{a\log(2)}{b}$ with $b\geq 1$ and 
	let $\Omega$ be a cycle of $T$, then there exists an integer  $n_0>0$ such that  for all $n\geq n_0$, we have
	$$  K\geq\min\Bigl(q_n,\dfrac{c_0 \min(\Omega_1)}{q_n^\mu}\Bigr),$$
	which implies,
	$$\#\Omega\geq\min\Bigl(q_n,\dfrac{c_0 \min(\Omega)}{q_n^\mu}\Bigr).$$
\end{theorem}
\begin{proof}
Let   $p_n/q_n$ be the convergents to $\xi$. Then, the sequence  $(q_n)$ is strictly increasing unbounded sequences. It follows that there exists an integer $n_0>0$ such that $q_{n_0}>\widehat{q}_0$ and  $\forall n>n_0$, we have $q_n>q_{n_0}>\widehat{q}_0$. So, according to (\ref{muequation}), we have
$\Bigl|p_n-q_n\xi\Bigr|>\dfrac{1}{q_n^{\mu}}$.
If $K\geq q_n$, the required inequality is trivial. Then, suppose $K<q_n$, it follows from Lemma~\ref{Lemma1ContinuedFractions}  that
$\Bigl|(K+L)-K\xi\Bigr|\geq \Bigl|p_n-q_n\xi\Bigr|$.  According to (\ref{log(a)log(2)}), we get
$$\dfrac{bK}{a\min(\Omega_1)\log(2)}>\dfrac{1}{q_n^{\mu}}. $$ Which gives the required result.
\end{proof}

\section{Bounds for $m$-circuits}
Following the terminology of Davidson, see \cite{Davidson1976, Lagarias1985} a cycle $$\Omega(x_0)=(x_0\longrightarrow T(x_0)\longrightarrow T^{(2)}(x_0)  \longrightarrow\ldots \longrightarrow T^{(k)}(x_0)=x_0   ),$$ of length $k$ is said to be a circuit  if there is a value $i$ for which
$$x_0>T(x_0)>\ldots>T^{(i)}(x_0)$$
and 
$$T^{(i)}(x_0)<T^{(i+1)}(x_0)<\ldots<T^{(k)}(x_0)=x_0.$$
For the classical sequence of Collatz,  Steiner \cite{Steiner1978} showed that  the only cycle that is a circuit is the trivial cycle.  Here, we have assumed that the sequence starts with an odd integer $x_0$, increases in $i+1$ steps until an even number is encountered, then the sequence decreases in $k-i+1$ steps until the odd number $x_0$ is again encountered.  In  such a situation, we will say that the cycle is of one oscillation. In \cite{SimonsWegner2005}, Simons and Weger  derived lower and upper bounds for hypothetical cycle with  $m$ oscillations.  A cycle with $m$ oscillations  contains $m$  local minimum odd integers $x_0$, $x_1,\ldots,x_{m-1}$  and $m$  local maximum even integers  $y_0$, $y_1,\ldots,y_{m-1}$ such that the sequence starts with the odd integer $x_0$, increases in $k_0$ steps until the even number $y_0$  is encountered, then the sequence decreases in $\ell_0$ steps, until the odd integer $x_1$ is encountered, again increases in $k_1$ steps until the even number $y_1$ is encountered, and so on until the odd integer $x_{m-1}$ is encountered. Then, the sequence increases in $k_{m-1}$ steps until the even integer $y_{m-1}$ is encountered and finally the sequence increases in $\ell_m$ steps until the integer $x_0$ is again encountered.   The number $K$ of odd integers  and the number $L$  of even numbers in the cycle $\Omega$ are given by
$$K=\sum_{i=0}^{m-1}k_i,\quad \textrm{and}\quad L=\sum_{i=0}^{m-1}\ell_i,.$$

The following Lemma, as we will see, is easy to obtain.
\begin{lemma}
	For $ i=0,\ldots, m-1$, we have
	\begin{equation}\label{m-circuitRelation}
		y_i=T^{(k_i)}(x_i)=\Bigl(\dfrac{a}{2}\Bigr)^{k_i}x_i+\dfrac{b}{a-2}\Bigl(\Bigl(\dfrac{a}{2}\Bigr)^{k_i}-1\Bigr)=2^{\ell_i}x_{i+1}, 
	\end{equation}
with the conditions  $x_m=x_0$, $k_m=k_0$ and $\ell_m=\ell_0$.
\end{lemma}
\begin{proof} 	For $ i=0,\ldots, m-1$, we have
	$T^{(j)}(x_i)$ is odd for $j=0,\ldots,k_i-1$. Then
	$T(x_i)=\dfrac{a}{2}x_i+\dfrac{b}{2}$ and
	$T^{(2)}(x_i)=\dfrac{a}{2}\left(\dfrac{a}{2}x_i+\dfrac{b}{2}\right)+\dfrac{b}{2}=
	\left(\dfrac{a}{2}\right)^2x_i+\dfrac{b}{2}\left(\dfrac{a}{2}+1 \right)$. Then,  by induction we obtain that, for $j=1,\ldots,k_i-1$, we have
	$$T^{(j)}(x_i)=\left(\dfrac{a}{2}\right)^jx_i+\dfrac{b}{2}\sum_{\ell=0}^{j-1}\left(\dfrac{a}{2} \right)^\ell,$$
	wish is also odd. Then, 
	$$y_i=T^{(k_i)}(x_i)=\dfrac{a}{2}\left(\left(\dfrac{a}{2}\right)^{k_i-1}x_i+\dfrac{b}{2}\sum_{\ell=0}^{k_i-2}\left(\dfrac{a}{2}
	\right)^{\ell}\right)+\dfrac{b}{2}.$$
	It follows that
	$$y_i=T^{(k_i)}(x_i)=\left(\dfrac{a}{2}\right)^{k_i}x_i+\dfrac{b}{2}\sum_{\ell=0}^{k_i-1}\left(\dfrac{a}{2}
\right)^{\ell}=\left(\dfrac{a}{2}\right)^{k_i}x_i+
\dfrac{b}{2}\left(\dfrac{\left(\dfrac{a}{2} \right)^{k_i}-1}{\dfrac{a}{2}-1}\right),$$
wish get the relation (\ref{m-circuitRelation}).
\end{proof}

\begin{lemma}\label{Lemma_moscilations}
	Let $b\geq 1$. A necessary condition  that a nontrivial cycle $\Omega$ having $m$ oscillations exists is 
	\begin{equation}\label{m-oscillations1}
		0<(K+L)-K\dfrac{\log(a)}{\log(2)}<\dfrac{b}{(a-2)\,\log(2)}\sum_{i=0}^{m-1}\dfrac{1}{x_i},
	\end{equation}
wish  implies that
	\begin{equation}\label{m-oscillations2}
	0<(K+L)-K\dfrac{\log(a)}{\log(2)}<\dfrac{m\,b}{(a-2)\,\log(2)\,\min(\Omega)}.
\end{equation}
\end{lemma}
\begin{proof}
	From the relation (\ref{m-circuitRelation}), we have
	$$\dfrac{2^{k_i+\ell_i}}{a^{k_i}}  \dfrac{x_{i+1}}{x_{i}}=1+\dfrac{b}{(a-2)x_i}\Bigl(1-\Bigl(\dfrac{2}{a}\Bigr)^{k_i}\Bigr), \quad i=0,\ldots, m-1.$$
	Taking the product from $i=0$ to $i=m-1$, and using the fact that $x_m=x_0$, we get
$$1<\dfrac{2^{K+L}}{a^{K}}  =\prod_{i=0}^{m-1}\left(1+\dfrac{b}{(a-2)x_i}\Bigl(1-\Bigl(\dfrac{2}{a}\Bigr)^{k_i}\Bigr)\right)<\prod_{i=0}^{m-1}\left(1+\dfrac{b}{(a-2)x_i}\right).$$
Now, applying the $\log$  and  using the inequality $\log(1+x)<x$ for $x>0$ we obtain the relation (\ref{m-oscillations1}) which implies obviously  
the relation (\ref{m-oscillations2}).
\end{proof}

\begin{theorem}
	Let us set   $c_1 =\dfrac{(a-2)\log(2)}{b}$ with $b\geq 1$. If  a nontrivial cycle $\Omega$ with $m$ oscillations exists, then for all $n\geq 1$, 
	$$ m\geq\min\Bigl(q_n,\dfrac{c_1 \min(\Omega)}{q_n+q_{n+1}}\Bigr).$$
\end{theorem}
\begin{proof}
	It is similar to the proof of Theorem~\ref{theorem1_K>=} by using Lemma~\ref{Lemma_moscilations}.
\end{proof}

\begin{theorem}
	Let us set   $c_1 =\dfrac{a\log(2)}{b}$ with $b\geq 1$.  If  a nontrivial cycle $\Omega$ with $m$ oscillations exists, then for all $n\geq 1$,  there exists an integer  $n_0>0$ such that  for all $n\geq n_0$, we have
	$$m\geq\min\Bigl(q_n,\dfrac{c_1 \min(\Omega)}{q_n^\mu}\Bigr),$$
	where $\mu$ is a measure of irrationality of $\dfrac{\log(a)}{\log(2)}$.
\end{theorem}
\begin{proof}
	It is similar to the proof of Theorem~\ref{theorem2_K>=} by using Lemma~\ref{Lemma_moscilations}.
\end{proof}

\section{Particular Syracuse problems}
\subsection{The $(2^\nu+1)n+(2^\nu-1)$ problem}
For the choice of the parameters  $a=2^\nu+1$ and $b=2^\nu-1$ for  $\nu\geq 1$. The corresponding map $T$ is  given by
$$
T(n)=\left\{\begin{array}{lc}
	n/2&\mathrm{if }\; $n$\; \mathrm{ is\; even},\\
	((2^\nu+1) n+(2^\nu-1) )/2&\mathrm{if} \;$n$\;\mathrm{ is\;odd}.\\
\end{array}\right.
$$
We have $a+b=2^{\nu+1}$, then according to Theorem~\ref{SetOfCycles}, the set $\mathcal{C}(2^\nu+1,2^\nu-1)$  contains at least   the  trivial cycle $\Omega(1)$ of length $\nu+1$, where
$$\Omega(1)=(1\rightarrow 2^{\nu}\rightarrow 2^{\nu-1} \rightarrow\ldots\rightarrow 2\rightarrow1).$$

For $\nu=1$ we recover  the classical Syracuse problem with $a=3$ and $b=1$ and we have the classical  Collatz conjecture, see Conjecture~\ref{Collatzconjecture}.

\begin{conjecture}Let  $\nu=2$ and  let the parameters   $a=2^\nu+1=5$ and $b=2^\nu-1=3$, then     the set  $\mathcal{C}(5,3)$ contains exactly $7$ cycles, we have $$\mathcal{C}(5,3)=\left\{\Omega(1),\Omega(3),\Omega(39),\Omega(43),\Omega(51),\Omega(53),\Omega(61)\right\},$$
and 
$$\mathbb{N}=G(1)\cup G(3)\cup G(39)\cup G(43)\cup G(51)\cup G(53)\cup G(61)\cup  G(\infty),$$
with $G(\infty)\not=\emptyset$.  The trivial cycle of length $3$ is
	$\Omega(1)=(1\rightarrow 4\rightarrow 2 \rightarrow 1)$, 
	the cycle $\Omega(3)$ of length $5$ is  $\Omega(3)=(3\rightarrow 9\rightarrow 24 \rightarrow 12\rightarrow  6 \rightarrow 3)$ and the other  $5$ nontrivial  cycles of length $7$  are:
	
	$\Omega(39)=(39\rightarrow99\rightarrow249\rightarrow624\rightarrow312\rightarrow156\rightarrow78\rightarrow39)$,
	
	$\Omega(43)=(43\rightarrow109\rightarrow274\rightarrow137\rightarrow344\rightarrow172\rightarrow86\rightarrow43)$,
	
	$\Omega(51)=(51\rightarrow129\rightarrow324\rightarrow162\rightarrow81\rightarrow204\rightarrow102\rightarrow51)$,
	
	$\Omega(53)=(53\rightarrow134\rightarrow67\rightarrow169\rightarrow424\rightarrow212\rightarrow106\rightarrow53)$
	and 
	
	$\Omega(61)=( 61\rightarrow154 \rightarrow77\rightarrow 194\rightarrow 97\rightarrow 244\rightarrow 122\rightarrow 61)$.
\end{conjecture}

\begin{conjecture}Let  $\nu\geq 3$ and let   the parameters  $a=2^\nu+1$ and $b=2^\nu-1$, then  the set  $\mathcal{C}(a,b)$ contains exactly one cycle. We have	$\mathcal{C}(a,b)=\left\{\Omega(1)\right\}$ and 
	$$\mathbb{N}=G(1)\cup  G(\infty),$$
	with $G(\infty)\not=\emptyset$ and the trivial cycle, of length $\nu+1$, is
	$$\Omega(1)=(1\rightarrow 2^\nu\rightarrow 2^{\nu-1}\rightarrow\ldots \rightarrow 2 \rightarrow1).$$ 
\end{conjecture}

As in the particular case, the problem in general case lies in the irregular behavior of the successive iterates.  First, we recall the definition of the expansion factor $s(n)$ given in \cite{Lagarias1985} by:
$$s(n)=\frac{ \displaystyle \sup_{k\geq 0}T^{(k)}(n)}{n},$$
if $n$ has  a  bounded trajectory $\Gamma(n)$  and $s(n)=+\infty$ if $\Gamma(n)$ is a divergent trajectory. 
The following result gives an idea on the computational difficulties encountered also in the general problem.  It  is similar to the one given in \cite{Crandall1978}.

\begin{theorem}
	For $\nu\geq 1$ and the parameters  $a=2^\nu+1$ and  $b=2^\nu-1$,  the sequence $(s(n))_{n\in\mathbb{N}}$ is unbounded.
\end{theorem}
\begin{proof}
For $k\in\mathbb{N}$, le $n_k=2^k-1$. We have immediately
$$T(n_k)=(a n_k+b )/2=a 2^{k-1}-1, \;\textrm{for}\;\;k> 1.$$
And 
$$T^{(2)}(n_k)=a^2 2^{k-2}-1, \;\textrm{for}\;\;k> 2.$$
Thus by induction
$$T^{(j)}(n_k)=a^j 2^{k-j}-1, \;\textrm{for}\;\;k> j.$$
Then,  for $k>1$, the number
$T^{(k-1)}(n_k)=a^{k-1} 2-1$
belongs to the trajectory $\Gamma(n_k)$. It follows that for $k>1$:
$$s(n_k)=\frac{ \displaystyle \sup_{k\geq 0}T^{(k)}(n_k)}{n_k}\geq \frac{ a^{k-1} 2-1}{2^k-1}>\Bigl(\frac{ a}{2}\Bigr)^{k-1},$$
as $a=2^\nu+1>2$ for $\nu\geq 1$, the right-hand side goes to $\infty$ as $k\longrightarrow\infty$.
\end{proof}

The following result follows the ideas given by  Rozier in \cite{Rozier1990} for the classical Syracuse problem.

\begin{theorem} Let  $\nu\geq 1$ and let the parameters  $a=2^\nu+1$ and  $b=2^\nu-1$, a
 necessary condition  that a nontrivial cycle $\Omega$  for $T$  having one  oscillation exists  is that there exists an integer $n_0>0$ such that for $K>n_0$,
	we have
	$$2^K-1\leq \dfrac{K^\mu}{\log(2)},\quad L=1-K+\Bigl\lfloor\dfrac{\log(a)}{\log(2)}K\Bigr\rfloor\quad\textrm{and}\quad \dfrac{a^K-2^K}{2^{K+L}-a^K}\in\mathbb{N},$$
	where $\mu$ is a measure of irrationality of $\dfrac{\log(a)}{\log(2)}$ ,  $K$ and $L$ are  the number of odd and even  integers in the cycle $\Omega$, respectively.
\end{theorem}
\begin{proof} 
	Suppose that  such a   nontrivial cycle $\Omega$ exists with one oscillation $m=1$. From the relation  (\ref{m-circuitRelation}), we have
	\begin{equation}\label{Relation,m=1}	y_0=T^{(K)}(x_0)=\Bigl(\dfrac{a}{2}\Bigr)^{K}x_0+\dfrac{b}{a-2}\Bigl(\Bigl(\dfrac{a}{2}\Bigr)^{K}-1\Bigr)=2^{L}x_{0}.
	\end{equation}
From the last relation and as $b=a-2$, we obtain that
$$x_0=\dfrac{a^K-2^K}{2^{K+L}-a^K}\in\mathbb{N}.$$
and
	$$1< \dfrac{2^{K+L}}{a^K}=1+\dfrac{1}{x_0}\Bigl(1-\dfrac{2^{K}}{a^K}\Bigr)<1+\dfrac{1}{x_0}.$$
	Applying the $log$ to the last inequality and taking in account into the inequality $\log(1+x)<x$ for $x>0$, we get
\begin{equation}\label{Relation_ROZIER1}
	0<\dfrac{K+L}{K}-\dfrac{\log(a)}{\log(2)}<\dfrac{1}{x_0\log(2)}\dfrac{1}{K},
\end{equation}
which may be also written as
	$$\dfrac{\log(a)}{\log(2)}K<K+L<\dfrac{\log(a)}{\log(2)}K+\dfrac{1}{x_0\log(2)},$$
	From the last relation and as $\min(\Omega)>2$,  ($\Omega$ is assumed to be nontrivial cycle) we get also that
	$$K+L=1+\Bigl\lfloor\dfrac{\log(a)}{\log(2)}K\Bigr\rfloor.$$
The following relation, from the relation (\ref{Relation,m=1}), 
$$2^K(x_02^{L}+1)=a^K(x_0+1).$$
shows that $2^K$ divise $x_0+1$, then 
\begin{equation}\label{Relation_ROZIER2}
x_0>2^K-1.
\end{equation}

From the definition of the measue of rationality $\mu$  of  $\dfrac{\log(a)}{\log(2)}$, there exists an integer $n_0>0$ such that for all pair  $(p,q)\in\mathbb{Z}\times \mathbb{N}$ with $q>n_0$ we have
$$\left| \dfrac{p}{q}-\dfrac{\log(a)}{\log(2)}\right|>\dfrac{1}{q^{\mu+1}}.$$
Then for $K>n_0$, we have
$$\dfrac{K+L}{K}-\dfrac{\log(a)}{\log(2)}>\dfrac{1}{K^{\mu+1}}.$$
It follows from (\ref{Relation_ROZIER1}) that,
$$\dfrac{1}{x_0\log(2)}\dfrac{1}{K}>\dfrac{K+L}{K}-\dfrac{\log(a)}{\log(2)}>\dfrac{1}{K^{\mu+1}}.$$
which gives from (\ref{Relation_ROZIER2}) that
$$\dfrac{K^\mu}{\log(2)}>2^K-1.$$
This concludes the proof.
\end{proof}

The  following corollary  gives  the main result given in \cite{Rozier1990}.  
\begin{corollary}
	For $\nu=1$, we get the classical case of the Syracuse problem corresponding to $a=3$ and $b=1$. In this case, it follows that there is no cycle with only one oscillation other than the trivial cycle.
\end{corollary}
\begin{proof}
	From the transcendental number theory, see  \cite{Rhin1987,Waldschmidt1987}, we have
	$$\left|\dfrac{p}{q}-\dfrac{\log(3)}{\log(2)}\right|>\dfrac{1}{q^{15}},\quad \forall (p,q)\in\mathbb{N}^2, \;\textrm{with}\;  q\geq 2.$$
	Then $\mu+1=15$, which implies that $\mu=14$. According to the previous Theorem, it follows that for the classical Syracuse problem,  necessary condition is that for $K\geq 2$, we have
		$$2^K-1\leq \dfrac{K^{14}}{\log(2)},\quad L=1-K+\Bigl\lfloor\dfrac{\log(a)}{\log(2)}K\Bigr\rfloor\quad\textrm{and}\quad x_0=\dfrac{3^K-2^K}{2^{K+L}-3^K}\in\mathbb{N},$$
		The condition  $2^K-1\leq \dfrac{K^{14}}{\log(2)}$ holds if and  only if $K\leq 91$ and  there is no integer $K\leq 91$ for which  $x_0=\dfrac{3^K-2^K}{2^{K+L}-3^K}$ is an integer for $L=1-K+\Bigl\lfloor\dfrac{\log(a)}{\log(2)}K\Bigr\rfloor$. 		
\end{proof}

The following Lemme is obviously  obtained from the expression of the parameter $A_\nu$ given in the Lemma.

\begin{lemma} For all integer $\nu\geq 1$, we set 
	$\displaystyle A_\nu=\log\Bigl(\dfrac{(2^\nu+1)^\nu-1}{2^\nu-1}\Bigr)/\log(2)=\log\Bigl(\dfrac{a^\nu-1}{b}\Bigr)/\log(2).$
	Then 
	\begin{itemize}
		\item We have  $a^\nu-1=2^{A_\nu}b$.
		\item For $\nu =1$,  $A_\nu=\nu\times (\nu-1)+1=1$ and for $\nu =2$, $A_\nu=\nu\times (\nu-1)+1=3$.
		\item  For all integer $\nu\geq 3$, we have  $\lfloor A_\nu\rfloor = \nu(\nu-1)$. More precisely,
		$$\nu(\nu-1)<A_\nu<\nu(\nu-1)+1.$$
		\item As $\nu\longrightarrow\infty $, we have   $A_\nu\sim \nu(\nu-1)$.
		
	\end{itemize}

\end{lemma}

\begin{theorem}
	The cycle $\Omega(b)$ of length   $\nu+A_\nu$  where
	
	\begin{eqnarray*}
\Omega(b)=(b\rightarrow a2^{\nu-1}-1\rightarrow a^22^{\nu-2}-1\rightarrow\ldots\rightarrow\\ a^\nu-1=2^{A_\nu}b\rightarrow 2^{A_\nu-1}b\rightarrow 2^{A_\nu-2}b\rightarrow \ldots\rightarrow 2b\rightarrow b),
	\end{eqnarray*}
exists only for the both cases where $A_\nu$ is an integer, namely for   $\nu=1$ or $\nu=2$.
\end{theorem}
\begin{proof}
	We have $T^{(\nu)}(b)=a^\nu-1$. According to the previous Lemma $a^\nu-1=2^{A_\nu}b$, which gives in the case where $A_\nu$ is an integer the relation
		$T^{(\nu+A_\nu)}(b)=b$. But  $A_\nu$ is an integer if and only if $\nu=1$ or $\nu=2$. For $\nu=1$, we have $b=2^\nu-1=1$ and $\Omega(b)=\Omega(1)$ is the trivial cycle. For $\nu=2$, we have $b=3$,  $A_\nu=3$ and the cycle
			$\Omega(b)$ is of length $\nu+A_\nu=2+3=5$,  is exactly the one given in the theorem, we have
		\begin{eqnarray*}
		\Omega(b)=(3\rightarrow a2^{\nu-1}-1=9\rightarrow a^2-1=24=2^\nu b=2^3\times 3\\
		\rightarrow 2^{A_\nu-1}b=2^2\times 3\rightarrow 2b=2\times 3\rightarrow 3),
		\end{eqnarray*}
			$$\Omega(b)=(3\rightarrow 9\rightarrow 24\rightarrow 12\rightarrow 6\rightarrow 3).$$		
		
\end{proof}

\subsection{The $(2^\nu-1)n+(2^\nu+1)$ problem}
For the choice of the parameters  $a=2^\nu-1$ and $b=2^\nu-1$ for  $\nu\geq 1$. The corresponding map $T$ is  given by
$$
T(n)=\left\{\begin{array}{lc}
	n/2&\mathrm{if }\; $n$\; \mathrm{ is\; even},\\
	((2^\nu-1) n+(2^\nu+1) )/2&\mathrm{if} \;$n$\;\mathrm{ is\;odd}.\\
\end{array}\right.
$$
We have $a+b=2^{\nu+1}$, then according to Theorem~\ref{SetOfCycles}, the set $\mathcal{C}(2^\nu-1,2^\nu-1)$  contains at least two trivial cycles.
The  first one $\Omega(1)$ of length $\nu+1$, where
$$\Omega(1)=(1\rightarrow 2^{\nu}\rightarrow 2^{\nu-1}\rightarrow 2^{\nu-2}\rightarrow\ldots 2\rightarrow 1),$$ 
and the second one $\Omega(b)$ of length $\nu$, where
$$\Omega(b)=(b\rightarrow 2^{\nu-1}b\rightarrow 2^{\nu-2}b\rightarrow \ldots\rightarrow 2b\rightarrow b).$$

\begin{conjecture}For $\nu=2$,  $a=2^\nu-1=3$ and $b=2^\nu+1=5$, the set $\mathcal{C}(a,b)$ contains  exactly $6$ cycles
	$\mathcal{C}(a,b)=\{\Omega(1),\Omega(5),\Omega(19),\Omega(23),\Omega(187),\Omega(347)\}$, and
	$$
	\mathbb{N}=G(1)\cup G(5)\cup G(19)\cup  G(23)\cup G(187)\cup G(347).$$
We have $G(\infty)=\emptyset$ with the  cycles of length $3$, $2$, $5$, $5$, $27$ and $27$, respectively,  are
	$$\begin{array}{rcl}
		\Omega(1)&=&(1\rightarrow4\rightarrow2 \rightarrow 1),\\
		\Omega(5)&=&(5\rightarrow 10 \rightarrow 5),\\
		\Omega(19)&=&(19\rightarrow 31\rightarrow 49\rightarrow 76\rightarrow 38\rightarrow 19),\\
		\Omega(23)&=&(23\rightarrow  37\rightarrow  58 \rightarrow 29\rightarrow  46\rightarrow  23),\\
		\Omega(187)&=&(187\rightarrow 283\rightarrow 427\rightarrow 643\rightarrow 967\rightarrow 1453\rightarrow 2182\rightarrow 1091\rightarrow 1639\rightarrow\\ && 2461
		\rightarrow 3694\rightarrow 1847\rightarrow 2773\rightarrow 4162\rightarrow 2081\rightarrow 3124\rightarrow 1562\rightarrow 781\rightarrow\\&& 1174\rightarrow 587\rightarrow 883\rightarrow 1327\rightarrow 1993\rightarrow 2992\rightarrow 1496\rightarrow 748\rightarrow 374\rightarrow 187),\\
		\Omega(347)&=&(347\rightarrow 523\rightarrow 787\rightarrow 1183\rightarrow 1777\rightarrow 2668\rightarrow 1334\rightarrow 667\rightarrow 1003\rightarrow\\
		&& 1507\rightarrow 2263\rightarrow 3397\rightarrow 5098\rightarrow 2549\rightarrow 3826\rightarrow 1913\rightarrow 2872\rightarrow 1436\rightarrow\\
		&& 718\rightarrow 359\rightarrow 541\rightarrow 814\rightarrow 407\rightarrow 613\rightarrow 922\rightarrow 461\rightarrow 694\rightarrow 347).\\
	\end{array}$$
\end{conjecture}

\begin{conjecture}For  $\nu\geq 1$ with $\nu\not=2$, $a=2^\nu-1$ and $b=2^\nu+1$, the set  $\mathcal{C}(a,b)$contains exactly two cycles. We have $$
	\mathcal{C}(a,b)=\left\{\Omega(1),\Omega(b)\right\},\;\;
	\textrm{and}\;\;
	\mathbb{N}=G(1)\cup G(b) \cup   G(\infty),$$
	with $G(\infty)=\emptyset$ for $\nu=1$ and $G(\infty)\not=\emptyset$ for $\nu\geq 3$. The trivial cycle of length $\nu+1$ is
	$$\Omega(1)=(1\rightarrow 2^{\nu}\rightarrow 2^{\nu-1}\rightarrow 2^{\nu-2}\rightarrow\ldots 2\rightarrow 1),$$ 
	and  the second cycle of length $\nu$ is 
	$$\Omega(b)=(b\rightarrow 2^{\nu-1}b\rightarrow 2^{\nu-2}b\rightarrow \ldots\rightarrow 2b\rightarrow b).$$ 
	
\end{conjecture}

\subsection{The $(2^\nu+1)n-(2^\nu-1)$ problem}
For the choice of the parameters  $a=2^\nu+1$ and $b=-2^\nu+1<0$ for  $\nu\geq 1$. The corresponding map $T$ is  given by
$$
T(n)=\left\{\begin{array}{lc}
	n/2&\mathrm{if }\; $n$\; \mathrm{ is\; even},\\
	((2^\nu+1) n-(2^\nu-1) )/2&\mathrm{if} \;$n$\;\mathrm{ is\;odd}.\\
\end{array}\right.
$$
We have $a+b=2$, then according to Theorem~\ref{SetOfCycles}, the set $\mathcal{C}(2^\nu+1,-2^\nu+1)$  contains at least two trivial cycles.
The  first one $\Omega(1)$ of length $1$, where
$$\Omega(1)=(1\rightarrow1),$$ 
and the second one $\Omega(-b)$ of length $\nu$, where
$$\Omega(-b)=(-b\rightarrow-b2^{\nu-1}\rightarrow-b2^{\nu-2}\rightarrow \ldots\rightarrow-2b\rightarrow-b).$$

\begin{conjecture}Let  $\nu=1$ and let   the parameters  $a=2^\nu+1=3$ and $b=-2^\nu+1=-1$, then  the set $\mathcal{C}(a,b)$ contains  exactly three cycles
	$\mathcal{C}(a,b)=\{\Omega(1),\Omega(5),\Omega(17)\}$, and
	$$
	\mathbb{N}=G(1)\cup G(5)\cup G(17).$$
	We have $G(\infty)=\emptyset$ with the  cycles of length $1$, $3$ and $11$, respectively,  are
	$$\begin{array}{rcl}
		\Omega(1)&=&(1\rightarrow 1),\\
		\Omega(5)&=&(5\rightarrow 7\rightarrow 10 \rightarrow 5),\\
		\Omega(17)&=&(17\rightarrow 25\rightarrow 37\rightarrow 55\rightarrow 82\rightarrow 41 \rightarrow 61\rightarrow 91\rightarrow 136\rightarrow 68\rightarrow 34 \rightarrow 17).\\\\
	\end{array}$$
\end{conjecture}

\begin{conjecture}\label{ParticularConjecture-nu>=2-Bis}Let  $\nu\geq 2$ and let   the parameters  $a=2^\nu+1$ and $b=-2^\nu+1$, then  the set  $\mathcal{C}(a,b)$contains exactly two cycles. We have $$
	\mathcal{C}(a,b)=\left\{\Omega(1),\Omega(-b)\right\},\;\;
	\textrm{and}\;\;
	\mathbb{N}=G(1)\cup G(-b) \cup   G(\infty),$$
	with $G(\infty)\not=\emptyset$. The trivial cycle of length $1$ is
	$\Omega(1)=(1\rightarrow1)$ and  the second cycle of length $\nu$ is 
	$$\Omega(-b)=(-b\rightarrow-b2^{\nu-1}\rightarrow-b2^{\nu-2}\rightarrow \ldots\rightarrow-2b\rightarrow-b).$$ 
	
\end{conjecture}

\subsection{The $(2^\nu-1)n+1$ problem}
For the choice of the parameters  $a=2^\nu-1$ and $b=1$ for  $\nu\geq 2$. The corresponding map $T$ is  given by
$$
T(n)=\left\{\begin{array}{lc}
	n/2&\mathrm{if }\; $n$\; \mathrm{ is\; even},\\
	((2^\nu-1) n+1 )/2&\mathrm{if} \;$n$\;\mathrm{ is\;odd}.\\
\end{array}\right.
$$
We have $a+b=2^{\nu}$, then according to Theorem~\ref{SetOfCycles}, the set $\mathcal{C}(2^\nu-1,1)$  contains at least   the  trivial cycle $\Omega(1)$, of length $\nu$, where
$$\Omega(1)=(1\rightarrow 2^{\nu-1}\rightarrow 2^{\nu-2} \rightarrow\ldots\rightarrow 2\rightarrow1).$$

For $\nu=2$ we recover  the classical Syracuse problem with $a=3$ and $b=1$ and we have the classical  Collatz conjecture, see Conjecture~\ref{Collatzconjecture}.

\begin{conjecture}\label{ParticularConjecture-2-nu>=3_2power_nu-1}Let  $\nu\geq 3$ and let   the parameters  $a=2^\nu-1$ and $b=1$, then  the set  $\mathcal{C}(a,b)$ contains exactly one cycle, $\mathcal{C}(a,b)=\left\{\Omega(1)\right\}$. We have 
$$
	\mathbb{N}=G(1)\cup  G(\infty),$$
	with $G(\infty)\not=\emptyset$ where  the trivial cycle $\Omega(1)$,  of length $\nu$ is
	$$\Omega(1)=(1\rightarrow 2^{\nu-1}\rightarrow 2^{\nu-2}\rightarrow\ldots \rightarrow 2 \rightarrow1).$$
\end{conjecture}





\end{document}